\documentclass[12pt]{amsart}
\usepackage{graphicx}
\usepackage{graphics}
\usepackage{amssymb}
\usepackage{amsxtra}

\input diagrams

\usepackage[arrow, matrix]{xy}

\theoremstyle{plain}

\newtheorem{theo}{Theorem}[section]
\newtheorem{lemm}[theo]{Lemma}
\newtheorem{prop}[theo]{Proposition}
\newtheorem{coro}[theo]{Corollary}

\newtheorem{conj}[theo]{Conjecture}
\newtheorem{prob}[theo]{Problem}
\theoremstyle{definition}

\newtheorem{defi}[theo]{Definition}
\newtheorem{rema}[theo]{Remark}

\newfont{\rmm}{cmr10 scaled 1000}
\newfont{\itt}{cmsl10 scaled 1000}

\newfont{\rM}{cmr10 scaled 1700}

\newcommand{\lb}{\label}
\newcommand{\mlb}{\label}

\newcommand{\rrf}[1]{(\ref{#1})}
\newcommand{\mrf}[1]{\ref{#1}}

\newarrow{dashto}{}{dash}{}{dash}>

\newarrow{mapsto}|--->

\newarrow{dotsto}|...>

\newcommand{\g}{\gamma}

\newcommand{\ve}{\varepsilon}

\renewcommand{\l}{\lambda}

\renewcommand{\r}{\rho}

\newcommand{\D}{\Delta}

\renewcommand{\L}{\Lambda}

\newcommand{\MM}{{\mathcal M}}
\newcommand{\NN}{{\mathcal N}}

\newcommand{\CCC}{{\mathbf{C}} }

\newcommand{\FFF}{{\mathbf{F}} }

\newcommand{\RRR}{{\mathbf{R}}}

\newcommand{\ZZZ}{{\mathbf{Z}}}

\newcommand{\gC}{{\mathfrak{C}}}

\newcommand{\gK}{{\mathfrak{K}}}

\newcommand{\gT}{{\mathfrak{T}}}

\newcommand{\Wh}{\text{\rm Wh }}

\newcommand{\rk}{\text{\rm rk }}

\newcommand{\supp}{\text{\rm supp }}

\newcommand{\bere}{\begin{rema}}
\newcommand{\bede}{\begin{defi}}

\renewcommand{\beth}{\begin{theo}}
\newcommand{\bele}{\begin{lemm}}
\newcommand{\bepr}{\begin{prop}}
\newcommand{\beeq}{\begin{equation}}
\newcommand{\bega}{\begin{gather}}
\newcommand{\been}{\begin{enumerate}}

\newcommand{\bedee}{\begin{defii}}
\newcommand{\bethh}{\begin{theoo}}
\newcommand{\belee}{\begin{lemmm}}
\newcommand{\beprr}{\begin{propp}}

\newcommand{\beco}{\begin{coro}}

\newcommand{\beal}{\begin{aligned}}

\newcommand{\enre}{\end{rema}}

\newcommand{\enco}{\end{coro}}
\newcommand{\enpr}{\end{prop}}
\newcommand{\enth}{\end{theo}}
\newcommand{\enle}{\end{lemm}}
\newcommand{\enen}{\end{enumerate}}
\newcommand{\enga}{\end{gather}}
\newcommand{\eneq}{\end{equation}}
\newcommand{\enal}{\end{aligned}}

\newcommand{\bq}{\begin{equation}}
\newcommand{\bqq}{\begin{equation*}}

\renewcommand{\leq}{\leqslant}
\renewcommand{\geq}{\geqslant}

\newcommand{\bu}{\bullet}

\newcommand{\mxx}[1]{\quad\mbox{#1}\quad}
 
\newcommand{\wi}{\widetilde}

\newcommand{\ove}{\overline}

\newcommand{\wh}{\widehat}

\newcommand{\sm}{\setminus}

\newcommand{\sbs}{\subset}

\newcommand{\tens}[1]{\underset{#1}{\otimes}}

\newcommand{\starr}[1]{\underset{#1}{*}}

\newcommand{\Lxi}{{\wh \L}_\xi}

\newcommand{\wrt}{~with respect to}
\newcommand{\ho}{homomorphism}

\newcommand{\noconf}{~no~confusion~is~possible}

\newcommand{\Prf}{{\it Proof.\quad}}

\newcommand{\smo}{C^{\infty}}

\newcommand{\chart}{\Phi_p:U_p\to B^n(0,r_p)}
\newcommand{\atlas}{\{\Phi_p:U_p\to B^n(0,r_p)\}_{p\in S(f)}}

\newcommand{\pr}{\partial}

\newcommand{\qs}{\hfill\square}

\begin{document}

\title[Novikov homology and knots]
{Twisted Novikov homology and Circle-valued Morse theory for knots and links}
%({\it sixth draft\/})}
\author{Hiroshi Goda}
\address{Department of Mathematics,
Tokyo University of Agriculture and Technology,
2-24-16 Naka-cho, Koganei,
Tokyo 184-8588, Japan}
\email{goda@cc.tuat.ac.jp}
\author{Andrei V. Pajitnov}
%\date{\today}
\address{Laboratoire Math\'ematiques Jean Leray UMR 6629,
Universit\'e de Nantes,
Facult\'e des Sciences,
2, rue de la Houssini\`ere,
44072, Nantes, Cedex, France}
\email{pajitnov@math.univ-nantes.fr}

\thanks{
The first author is partially supported
by Grant-in-Aid for Scientific Research,
(No. 15740031),
Ministry of Education, Science, 
Sports and Technology, Japan.
}

\maketitle

\section{Introduction}
\mlb{s:intro}

Let $L$ be a link, that is, a $\smo$ embedding of disjoint union 
of the oriented circles in $S^3$.
Recall that  $L$ is called {\it fibred}
if there is a fibration $\phi:C_L=S^3\sm L\to S^1$
behaving ``nicely'' in a neighborhood of $L$.
If $L$ is not fibred, it is still possible to construct
a  Morse map $f:C_L=S^3\sm L \to S^1$
behaving nicely in a neighborhood of $L$;
such a  map has necessarily  a finite number of critical points.
The minimal number of  critical points 
of such map 
is an invariant of the link, called {\it Morse-Novikov number }
of $L$ and denoted by $\MM\NN(L)$;
it was first introduced 
and studied in \cite{prw}.
This invariant can be studied via the 
methods of the {\it Morse-Novikov theory}, in particular 
the Novikov inequalities \cite{novidok}  provide
the   lower estimates for the number $\MM\NN(L)$.
This inequalities can be considered as  {\it fibering obstructions}
for the link. 
The Alexander polynomials provide another fibering obstruction:
if a link is fibred then its Alexander polynomial is monic
(see \cite{rolfsen},  10.G.9).
The invariants arising from the Novikov theory
are in a sense stronger than the Alexander polynomial,
since they allow in certain cases to give 
lower bounds for $\MM\NN(L)$, see \cite{prw}
for examples of knots with arbitrarily large 
Morse-Novikov number.
Both the Alexander polynomial and the Novikov homology
above are {\it abelian invariants}, that is, 
they are calculated from the homology of the
 infinite cyclic covering of $C_L$.
More information is provided by the {\it non-abelian }
coverings, although the corresponding invariants are 
more difficult to deal with. 
Several non-abelian versions 
of the Alexander polynomial were developed in 90s
(see the papers by X.S.Lin \cite{lin}, M.Wada \cite{wada},
T.Kitano \cite{kitano}, and Kirk-Livingston \cite{kirklivingston}).
In the recent preprint \cite{gokimo}
by H.Goda, T.Kitano, and T.Morifuji, 
it was shown that if a knot is fibred, then
the twisted Alexander invariant is monic.

In the present paper we shall develop the 
twisted version of the Novikov homology.
This part (Section 3) may be of independent  interest for the Morse-Novikov theory.
The existing versions of the non-abelian Novikov homology for a CW complex $X$
are either modules over the Novikov completion of 
the $\ZZZ\pi_1(X)$, or the vector spaces associated 
to some representation of $\pi_1(X)$. Our version is in a sense 
intermediate between these two: the twisted Novikov homology
which we define is a module over the ring $\ZZZ((t))$, thus it is 
computable (as we shall show in this paper) 
and on the other hand 
%is contains the information about the torsion elements,
%which is important for geometric applications.
it allows to keep track of the non-abelian 
homological algebra associated to $\ZZZ\pi_1(X)$.
Theorem \mrf{t:nov_ineq} gives a lower bound for 
the Morse-Novikov number of a link in terms of the twisted Novikov
homology. 

We show that the twisted Novikov homology is additive
\wrt~ the connected sum of knots.
We apply these techniques to determine 
the Morse-Novikov numbers of the knots $n\gK\gT,n\gC$,
where  $\gK\gT$ is the Kinoshita-Terasaka knot \cite{kinoshitaterasaka}, 
$\gC$ is the Conway knot \cite{conway}, and $nK$ stands for the connected sum of
the knot $K$ with itself.   We show that
$$
\MM\NN(n\gK\gT) = \MM\NN(n\gC)\geq 2n/5.
$$
The computations of the invariants arising from the twisted Novikov homology 
for these knots was done with the help of 
the Kodama's KNOT program
(available at http://www.math.kobe-u.ac.jp/~kodama\\
/knot.html).
Applying the techniques of the papers \cite{godaone}, \cite{godatwo}
we prove that 
$$
\MM\NN(n\gK\gT) = \MM\NN(n\gC)\leq 2n.
$$
We recall the corresponding notions and 
results in Section 2 
from \cite{gabai}, \cite{godaone} and \cite{godatwo}.

We show that the twisted Novikov homology 
is related to twisted Alexander polynomials 
in the same way as the abelian Novikov homology is 
related to the classical Alexander polynomial.

The final section of the paper contains a discussion about
necessary and sufficient condition for a link to be fibred.

%We note that here 
%estimates of Morse-Novikov number by using several 
%geometric invariants are announced in \cite{hirarud}. 

Acknowledgements.
This paper was finished during the visit of 
the second author to Osaka City University. 
He is grateful to Professor Akio Kawauchi 
and Osaka City University for warm hospitality. 
Authors are grateful to Professor Teruaki Kitano, 
Professor Kouji Kodama and Professor Sadayoshi Kojima 
for valuable discussions.

%%%%%%%%%%%%%%%%%%%%%%%%%%%%%%%%%%%%%%%%%%%%%%%%%

\section{Heegaard splitting for sutured manifold} 
\mlb{s:suture}

The basic Morse theory gives a relationship 
of a Morse map and a handle decomposition for 
a manifold. 
In this section, we review the 
notion of Heegaard splitting for sutured manifold
introduced in \cite{godaone} and \cite{godatwo}, 
and reveal the relationship with circle-valued Morse maps. 
Moreover, we present some properties which are used to 
determine Morse-Novikov numbers. 

In this section, we assume that a link is always non-split.

Sutured manifold was defined by D.Gabai \cite{gabai}. 
We use here a special class of sutured manifolds, 
namely, complementary sutured manifolds:
 
\bede\mlb{d:suture}
Let $L$ be an oriented link in $S^3$, and 
$R$ a Seifert surface of $L$. 
Set $R_E=R\,\cap E(L)$ ($E(L)=$ cl$(S^3-N(L))$), 
and $(P, \delta)=(N(R_E,E(L)), N(\pr R_E, \pr E(L))$. 
We call $(P, \delta)$ a {\it product sutured manifold\/} 
for $R$. 
Thus $P$ is homeomorphic to $R_E\times [0,1]$, 
and then we denote by $R_{+}(\delta)$ ($R_{-}(\delta)$ resp.) 
the surface $R_E\times\{ 1\}$ ($R_E\times\{ 0\}$ resp.).  
Let 
$$(M, \g)=({\rm cl}(E(L)-P), {\rm cl}(\pr E(L)-\delta))$$ 
with $R_{\pm}(\gamma)=R_{\mp}(\delta)$.
We call $(M,\g)$ a 
{\it complementary sutured manifold\/} for $R$. 
In this paper, we call this a {\it sutured manifold\/} 
for short. 
\end{defi}

Here we denote by $N(X,Y)$ a regular neighborhood 
of $X$ in $Y$. 

In \cite{CG}, 
the notion of compression body was intorduced
by A. Casson and C. Gordon. 
It is a generalization of a handlebody, 
and important to define a Heegaard splitting 
for 3-manifolds with boundaries. 
 
\bede\mlb{d:compression}
A {\it compression body\/} $W$ is a cobordism 
rel $\pr$ between surfaces $\pr_{+}W$ and 
$\pr_{-}W$ such that 
$W\cong\pr_{+}W\times [0,1]\,\cup$ 2-handles $\cup$ 3-handles 
and $\pr_{-}W$ has no 2-sphere components. 
We can see that if $\pr_{-}W\neq\emptyset$ and $W$ is connected, 
$W$ is obtained from $\pr_{-}W\times [0,1]$ by attaching 
a number of 1-handles along the disks 
on $\pr_{-}W\times\{ 1\}$ 
where $\pr_{-}W$ corresponds to 
$\pr_{-}W\times\{ 0\}$. 

We denote the number of these 1-handles by $h(W)$. 
\end{defi}

These notions enable us to define a Heegaard splitting 
for sutured manifold. 

\bede\mlb{d:heegaard}(\cite{godaone})
$(W,W')$ is a Heegaard splitting for $(M,\g)$ if
\begin{enumerate}
\renewcommand{\labelenumi}{(\roman{enumi})}
\item 
$W,\,W'$ are compression bodies, 
\item
$W\cup W'=M$, 
\item
$W\cap W'=\pr_{+}W=\pr_{+}W', 
\pr_{-}W=R_{+}(\g),$ and $\pr_{-}W'=R_{-}(\g)$. 
\end{enumerate}
\end{defi}

\bede\mlb{d:genus}(\cite{godatwo})
Set 
$h(R)=\min\{h(W)(=h(W'))| (W,W')$ is a Heegaard splitting for the sutured manifold 
of $R\}$. 
We call $h(R)$ the {\it handle number\/} of $R$. 
%Further, we denote by $h(L)$ the minimum handle number 
%among all Seifert surfaces of $L$. 
\end{defi}

%Note that $L$ is a fibred if and only if $h(L)=0$.
The behaviors of the handle number are studied in 
\cite{godaone} and \cite{godatwo}.

In order to state the relationship between 
the handle number and the Morse-Novikov number, 
we recall some definitions on circle-valued Morse map 
according to \cite{prw}.

\bede\mlb{d:reg_morse}(\cite{prw})
Let $L$ be a link.
A Morse map $f:C_L\to S^1$ is said to be {\it regular}
if the link $L$
has a neighborhood framed as $L\times D^2$
with $L\approx L\times 0$,
in such a way that
the restriction $f|\: L\times (D^2 \sm \{0\})\to S^1$
is given by $(x,y)\mapsto y/|y|$.
\end{defi}

%Observe that for a regular Morse map $f:C_K\to S^1$
%we have necessarily
%$$f_*=\xi:G=\pi_1(C_K)\to\ZZZ=\pi_1(S^1).$$

For a regular 
Morse map  we denote by $S_i(f)$
the set of all critical points of index $i$
and by 
$m_i(f)$ the cardinality of $S_i(f)$.
We say that a Morse map 
$f: C_L\to S^1$ is {\it minimal\/} 
if it is regular and 
$m_i(f)$ are minimal possible among all 
regular maps homotopic to $f$.  
We define $\MM\NN(L)$ as the 
number of critical points of 
the minimal Morse map.

\bede\mlb{d:moderate}
A regular Morse map 
$f: C_L\rightarrow S^1$ is said to 
be {\it moderate\/} if it satisfies the all 
of the following:
\begin{enumerate}
\renewcommand{\labelenumi}{(\roman{enumi})}
\item 
$m_{0}(f)=m_{3}(f)=0$;
\item
all critical values corresponding to 
critical points of the same index coincide; 
\item
$f^{-1}(x)$ is a connected Seifert surface for 
any regular value $x\in S^{1}$. 
\end{enumerate}
\end{defi}

\beth\mlb{t:moderate}$($\cite{prw}$)$
Every link has a minimal Morse map
which is moderate. 
\enth

\beco
\begin{enumerate}
\renewcommand{\labelenumi}{(\theenumi)}
\item 
Let $f$ be a moderate map, then 
$m_1(f)=m_2(f)$. 
\item
Let $f$ be a regular Morse map 
realizing $\MM\NN(L)$, then 
$\MM\NN(L)=m_1(f)+m_2(f)$. 
\item
$\MM\NN(L)=2\times\min\{h(R)|R$ is a Seifert surface 
           for $L\}$. 
\end{enumerate}
\enco      

We denote by $h(L)$ the minimum handle number 
among all Seifert surfaces of $L$. 
Note that $L$ is a fibred if and only if $h(L)=0$.

Thus we know that the handle number and Morse-Novikov 
number are same essentially, that is, 
$$\MM\NN(L)=2\times h(L).$$ 

%From these observations, we know Conjecture 6.3 in
%\cite{prw} is true because of Theorem A in 
%\cite{godatwo} under the assumption that 
%the links are non-split.  
%On the other hand, 
%the question: 
%$\MM\NN(K_1\sharp K_2)=\MM\NN(K_1)+\MM\NN(K_2)?$ 
%is still open where $\sharp$ means 
%a connected sum. 
%Note that the inequality:
%\bq\lb{f:connect}
%\MM\NN(K_1\sharp K_2) \le \MM\NN(K_1)+\MM\NN(K_2)
%\end{equation}
%has been proved. 
%If $\MM\NN(L)$ is always attained by an 
%incompressible Seifert surface of $L$, then 
%we have the above equation. 

We shall finish this section with some remarks 
on behaviors of the invariant introduced above 
with respect to connected sum and plumbing. 
Let us denote $\natural$
the operation of plumbing. 
For a Seifert surface $R$ of a link $L$, 
we set $\MM\NN(R)=2\times h(R)$. 
Conjecture 6.3 in \cite{prw} says that 
$$\MM\NN(R_1\natural R_2)\le \MM\NN(R_1)+\MM\NN(R_2).$$
This conjecture follows from Theorem A 
in \cite{godatwo} in the case of non-split links. 
In the recent paper \cite{hirarud}, they 
prove it in general case. 
Let $K_1$ and $K_2$ be a knot in $S^3$, and 
we denote by $K_1\sharp K_2$ their connected sum.
As far as we know the question in \cite{prw} 
$$\MM\NN(K_1\sharp K_1)=\MM\NN(K_1)+\MM\NN(K_2)?$$ 
is still unproved, so we know only
\bq\lb{f:connect}
\MM\NN(K_1\sharp K_2) \le \MM\NN(K_1)+\MM\NN(K_2).
\end{equation}

\section{Twisted Novikov homology}
\mlb{s:twist_nov}

Let $R$ be a commutative ring, put 
$$
Q=R[t, t^{-1}],
\quad
\wh Q= R((t))=R[[t]][t^{-1}].
$$
The ring $Q$  is isomorphic to 
the group ring $R[\ZZZ]$, via the isomorphism
sending $t\in Q$ to the element $-1\in \ZZZ$.
The ring  
$\wh Q$ is then identified with 
{\it the Novikov completion } of 
$R[\ZZZ]$.
In the case when $R$ is a field, 
$\wh Q$ is also a field. 
When $R=\ZZZ$, the ring $\wh Q$ is PID. 
We shall need in this paper only the particular cases
when $R=\ZZZ$, or $R$ is a field. 
Let $X$ be a CW complex; let $G=\pi_1X$, and let
$\xi:G\to\ZZZ$ be a homomorphism.
Let $\rho:G\to GL(n,R)$
be a 
right 
representation of $G$, 
that is, 
$\r(g_1g_2)=\r(g_2)\r(g_1)$. 
The homomorphism $\xi$ extends to a ring homomorphism
$\ZZZ[G]\to Q$, which will be denoted by the same symbol $\xi$.
The tensor product $\r\otimes \xi$
(where $\xi$ is considered as a representation
$G\to GL(1, Q)$)
induces a representation
$\rho_\xi:G\to GL(n,Q)$.
The composition of this representation with
the natural inclusion
$Q\rInto \wh Q$
gives a representation
$
\wh\rho_\xi:G\to GL(n,\wh Q).
$
Let us  form a  chain complex
\bq\lb{f:chain}
\wh C_*(\wi{X};\xi,\r )=
\wh Q^n\tens{\wh\r_\xi} ~C_*(\wi{X}).
\end{equation}
Here $\wi{X}$ is the universal 
cover of $X$, 
$C_*(\wi{X})$ is a module over $\ZZZ\pi_1X$, 
and $\wh Q^n$ is a $\ZZZ G$-module 
via the representation $\wh\rho$.
Then (\ref{f:chain}) is a chain complex of free left 
modules over 
$\wh Q$, and the same is true for its
homology.
The modules
$$
\wh H_*({X};\xi,\r )=
H_*(\wh C_*(\wi{X};\xi,\r )),
$$
will be called
{\it $\rho$-twisted Novikov homology}
or simply 
 {\it twisted Novikov homology } if \noconf.
When these modules are finitely generated
(this is the case for example for
any $X$ homotopy equivalent to a finite CW complex)
we set
$$
\wh b_i(X; \xi,\r )= \rk_{\wh Q} (\wh H_i(X ;\xi,\r )),
\quad
\wh q_i(X ;\xi,\r )= {\rm t.n.}_{\wh Q} (\wh H_i(X ;\xi,\r )).
$$
where t.n. stands for the {\it torsion number }
of the $\wh Q$-module, 
that is the minimal possible number of generators
over $\wh Q$.

The numbers $\wh b_i(X;\xi,\r )$ and 
$\wh q_i(X;\xi,\r )$ can be 
recovered from the canonical decomposition of 
$\wh H_i(X;\xi,\r )$ into a direct sum 
of cyclic modules. 
Namely, let 
$$\wh H_i(X;\xi,\r )=\big(\bigoplus_{j=1}^{\alpha_{i}} \wh Q\big)
\bigoplus\big(\bigoplus_{j=1}^{\beta_{i}}\wh Q/\lambda^{(i)}_{j}\wh Q
\big)$$
where $\lambda^{(i)}_{j}$ 
are non-zero non-invertible elements of $\wh Q$ and 
$\lambda^{(i)}_{j+1}|\lambda^{(i)}_{j}\,\,\,\forall j$.
(Such decomposition exists since $\wh Q$ is a PID.)
Then $\alpha_{i}=\wh b_i(X;\xi,\r )$ and 
$\beta_{i}=\wh q_i(X;\xi,\r )$. 
It is not difficult to show that 
we can always choose $\lambda^{(i)}_{j}\in Q\,\,\,
\forall i,\,\forall j$. 

When $\r$ is a trivial 1-dimensional representation,
we obtain the usual Novikov homology,
which can be also calculated from the infinite cyclic covering
$\bar X$ 
associated to $\xi$,
namely
$$
\wh H_*(X;\xi,\r )
=
\wh Q\tens{Q} ~H_*(\bar X)
\mxx{ for }
\r=1:G\to GL(1, R).
$$

If $R$ is a field the numbers 
$\wh q_i(X;\xi,\r )$ vanish (for every representation $\rho$),
and the module $\wh H_*(X;\xi,\r )$
is a vector space over the field $\wh Q$.

%%%%%%%%%%%%%%%%%%%%%%%%%%%%%%%%%%%%%%%%%%%%%%%%

%\input nov_ineq.doc

\section{Novikov-type inequalities for knots and links}
\mlb{s:nov_ineq}

Now we shall apply the algebraic techniques 
developed in the previous section to
the topology of knots and links.
Let $L\sbs S^3$ be an oriented link and 
put $C_L=S^3\sm L$. Let $G$ denote $\pi_1(C_L)$.
%Let $\mu\in H_1(C_L)$ be the image of the positively oriented
%meridian of $K$; then 
For every positively oriented meridian $\mu_i$ 
of a component of $L$, 
%$\mu$ is the generator of the group
%$H_1(C_K)\approx \ZZZ$. 
there is a unique element
$\xi\in H^1(M, \ZZZ)$ such that 
for each $i$ we have $\xi(\mu_i)=1$.
We shall identify this element with the 
corresponding homomorphism
$G\to\ZZZ$.

%\bede\mlb{d:reg_morse}(\cite{prw})
%A Morse map $f:C_K\to S^1$ is said to be {\it regular}
%if the knot $K$
%has a neighborhood framed as $S^1\times D^2$
%with $K\approx S^1\times 0$,
%in such a way that
%the restriction $f|\: S^1\times (D^2\sm \{0\})\to S^1$
%is given by $(x,y)\mapsto y/|y|$.
%\end{defi}

%Observe that for a regular Morse map $f:C_K\to S^1$
%we have necessarily
%$$f_*=\xi:G=\pi_1(C_K)\to\ZZZ=\pi_1(S^1).$$

%For a regular 
%Morse map  we denote by $S_k(f)$
%the set of all critical points of index $k$
%and by 
%$m_k(f)$ the cardinality of $S_k(f)$.

Let $\r:G\to GL(n,R)$
be any representation of $G$ (where  $R=\ZZZ$ or $R$ is a field).
The next theorem follows from 
the main theorem in \cite{patou}. 
See also Theorem \mrf{t:nov_fibr} of the present paper. 
%an exercise in the Morse-Novikov theory
%(See \cite{patou} for the proof of a similar result).
\beth\mlb{t:nov_complex}
There is a free chain complex $\NN_*$ over $\wh Q$
such that 
\been\item
for every $i$ the number of 
free generators of $\NN_*$
in degree $k$ equals $n\times m_i(f)$;
\item
$H_*(\NN_*)\approx \wh H_*(C_L;\xi, \r)$.\enen
\enth

In what follows, we use the following terminology.
The cohomology class $\xi$ is
 determined by the orientation of the link, so
we shall omit it in the notation.
Further, 
$$\wh H_*(L,\r)=\wh H_*(C_L,\r),$$
$$\wh b_i(L,\r)=\wh b_i(C_L,\r),$$
$$\wh q_i(L,\r)=\wh q_i(C_L,\r).$$

The next theorem follows from Theorem \mrf{t:nov_complex}
by a simple algebraic argument.
\beth\mlb{t:nov_ineq}
Let $f:C_L\to S^1$
be any regular map.
Then 
\begin{multline}
m_i(f)
\geq 
\frac 1n\big(
\wh b_i(L,\r) +
\wh q_i(L,\r)+
\wh q_{i-1}(L,\r)
\big)\\
\mxx{ for every }
i.\hspace{1cm} \square
\end{multline}
\enth

\beco\mlb{c:obstr_fiber}
If $L$ is fibred, then $\wh H_*(L,\rho)=0$, and $\wh b_i(L, \r)=\wh q_i(L, \r)=0$
for every representation $\rho$ and every $i$.
\enco

\bepr\mlb{p:numbers}
The twisted Novikov numbers of the space $C_L$
satisfy the following relations:
\begin{gather}
\wh b_i(L, \r)
= \wh  q_i(L, \r)=\wh  q_2(L, \r)=0
\mxx{ for} i=0, i\geq 3,\quad \lb{f:numbers_1}\\
\wh b_1(L, \r)=\wh b_2(L, \r).\quad \lb{f:numbers_2}
\end{gather}
\enpr
\Prf
According to Theorem 3.3 of \cite{prw}
there is a regular map $f:C_L\to S^1$
such that $f$ has only critical points of indices 1 and 2
and $m_1(f)=m_2(f)$.
Using Theorem \mrf{t:nov_ineq}
we deduce \rrf{f:numbers_1}.
As for the point \rrf{f:numbers_2}
it follows from the fact that the Euler characteristics
of the chain complex $\NN_*$ is equal to $0$. $\qs$

In view of the preceding theorem the non-trivial part of the
Novikov inequalities is as follows:
\begin{gather}
m_1(f)\geq
\frac 1n \big(
 \wh b_1(L, \r)+ \wh q_1(L, \r)
\big)
;\\
m_2(f)\geq 
\frac 1n \big(
\wh b_1(L, \r)+ \wh q_1(L, \r)\big).
\end{gather}

Let us consider some examples and particular cases.

{\bf 1. } $R=\ZZZ$ and 
 $\rho$ is the  trivial 1-dimensional 
representation. 
The chain complex $\wh C_*(\wi C_L;\xi, \r)$
is equal to the chain complex
$C_*(\bar C_L)$,
where  $\bar C_L$ is the infinite cyclic covering associated to 
the cohomology class $\xi$.
Thus the twisted Novikov homology in this case coincides with the 
Novikov homology for links introduced and studied 
in \cite{prw}.
Theorem \mrf{t:nov_ineq}
and Proposition \mrf{p:numbers}
in this case are reduced to Proposition 2.1 
and the formulas (2) -- (6) of \cite{prw}.

{\bf 2. } The ring $R$  is a field.
In this case the Novikov ring $R((t))$ is also a field, and the torsion numbers
$\wh q_i(L,\r)$ vanish for every $i$ and every representation $\rho$.
The Novikov inequalities have the simplest possible form:
$$m_1(f)\geq \frac 1n 
\wh b_1(L, \r)\leq m_2(f).$$

{\bf 3.} 
Now we shall investigate the twisted Novikov homology 
for the connected sum of knots.
Let $K_1, K_2$ be oriented knots in $S^3$, and
put $K=K_1\sharp K_2$.
We have:
$$
\pi_1(K)= \pi_1(K_1) \starr{{Z}} \pi_1(K_2),
$$
where $
 {Z}$ is the infinite cyclic group generated by a meridian $\mu$ of
$K$ (see \cite{bz}, Ch.7, Prop. 7.10).
In particular
 the groups $\pi_1K_1, \pi_1K_2$ are naturally embedded into 
$\pi_1K$, and the meridian element  $\mu\in \pi_1K$ is the image of  some meridian
elements
$\mu_1, \mu_2$ of $K_1$, resp. $K_2$.
Now let $\rho_1:\pi_1K_1\to
GL(m,R),~\rho_2:\pi_1K_2\to GL(m,R)$
be two representations.
%Replacing $\rho_2$ by a conjugate representation, we can always consider
Assume 
that $\rho_1(\mu_1)= \rho_2(\mu_2)$.
In this case we obtain the product representation
$\rho=\r_1*\r_2:\pi_1K\to GL(m,R)$.
\beth\mlb{t:conn_sum}
$
\wh H_*(K,\r)
\approx
\wh H_*(K_1,\r_1)
\oplus
\wh H_*(K_2,\r_2).
$
\enth
\Prf
The complement $C_K$ is the unions of two subspaces
$C_1, C_2$ with $C_i$ having the homotopy type of $C_{K_i}$
(for $i=1, 2$).
The intersection $C_1\cap C_2$ is homeomorphic to the
twice punctured sphere $\D'=S^2\sm \{*,*\}$.
The universal covering
is therefore the union of two subspaces, which have the Novikov homology respectively
equal to $\wh H_*(K_1,\r_1)$
and $\wh H_*(K_2,\r_2)$.
The intersection of these two subspaces have the same Novikov homology as
$\D'$, and this module vanishes. Then the standard
application of the Mayer-Vietoris sequence
proves  the result sought. $\qs$

\beco\mlb{c:n_times}
Denote by $nK$ the connected sum of $n$ copies of the knot $K$.
Let $\r:\pi_1K\to GL(m,R)$ be a
representation (where $R=\ZZZ$ or $R$ is a field).
Let $\r^n:\pi_1(nK)\to GL(m,R)$ be the product of $n$ 
copies of representations $\r$.Then
$$\wh q_1(nK,\rho^n)=n\cdot \wh q_1(K,\rho).$$
\enco
 \Prf
This follows from the purely algebraic equality:
$$
{\rm t.n.}(nN)= n\cdot \big({\rm t.n.}(N)\big)
$$
where $N$ is any 
finitely generated module over a principal ideal domain, 
and $nN$ stands for the direct sum of $n$ copies of $N$.
$\qs$

%\input conway.doc

%%%%%%%%%%%%%%%%%%%%%%%%%%%%%%%%%%%%%%%%%%%%%%%%%%%%%

\section{The Kinoshita-Terasaka and Conway knots:
lower bounds for the Morse-Novikov numbers}
\mlb{s:kt_c}

The Kinoshita-Terasaka knot $\gK\gT$
was introduced in the paper \cite{kinoshitaterasaka},
and the Conway knot $\gC$ was discovered by J.Conway
much later \cite{conway}.
The resemblance of these two knots is
obvious on the pictures as in Figures 1 $\&$ 2
and is proved by coincidence
of many algebraic invariants
(still these knots are different,
as was proved by M.Wada in \cite{wada}).
\vskip0.2in

%$$
%\includegraphics{KTpic/ess.1}
%$$\vskip0.2in
%\centerline{ The Kinoshita-Terasaka knot $\gK\gT$}
%\vskip0.2in
%\vskip0.2in
%$$
%\includegraphics{conwpic/ess.1}
%$$\vskip0.2in
%\centerline{ The Conway knot $\gC$}
%\vskip0.2in

\begin{figure}
\centering
\includegraphics[width=.5\textwidth]{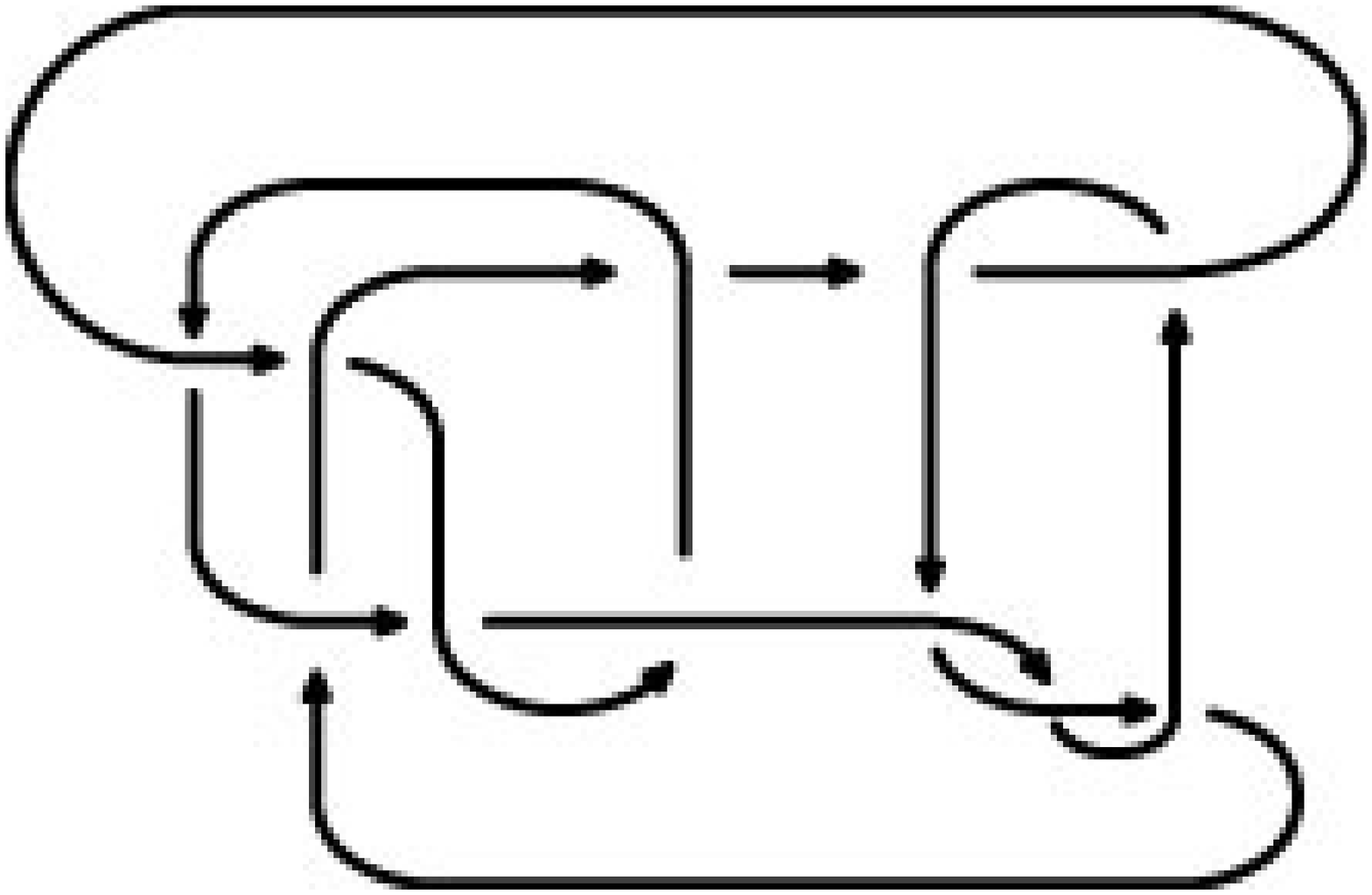}
\caption{The Kinoshita-Terasaka knot $\gK\gT$}
\end{figure}

\begin{figure}
\centering
\includegraphics[width=.5\textwidth]{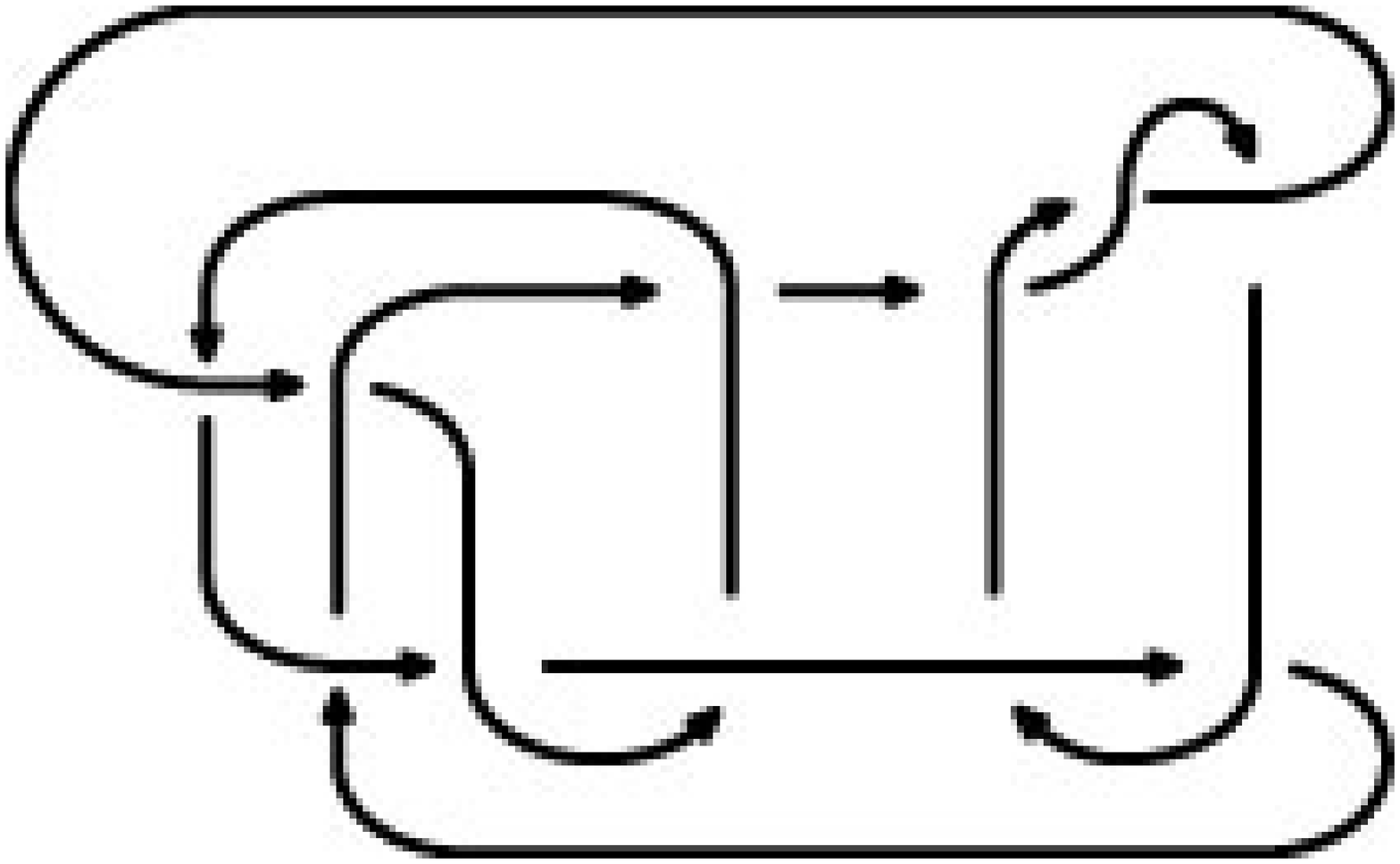}
\caption{The Conway knot $\gC$}
\end{figure}

%It was proved in \cite{gokimo}
%that the  Kinoshita-Terasaka knot
%is not fibred. 
%In this section we shall prove  that the Conway knot is not fibred, and
%moreover we shall prove that there is a representation
%$\rho:\pi_1\gC\to SL(5,\ZZZ)$
%such that
%$\wh q_1(\gC, \rho)\not=0.$
%(Note that we can also know that these knots are not fibred 
%from the known equality concerning the degree of
%Alexander polynomial and the genus of a fibred knot.)
%This implies by the Corollary \mrf{c:n_times}
%that  for every $n\geq 1$ we have

These knots are not fibred. 
Indeed, for a fibred knot the degree of its 
Alexander polynomial equals to two times of 
the genus of the knot, and the Alexander polynomial 
is trivial for both knots. 

In this section, we prove the following: 

\beth\mlb{t:conway}
There is a representation 
$\rho:
\pi_1\gC\to SL(5, \ZZZ)$ 
such that $\wh q_1(\gC,\r)\not= 0$.
\enth

By Corollary \mrf{c:n_times}, 
this theorem implies

\bq\lb{f:conw_estimate}
\MM\NN(n\gC) \geq 2n/5
\mxx{ for every } n.
\end{equation}

\Prf
The Wirtinger presentation for the group $\pi_1\gC$
has 11 generators $s_i$
and 11 relations:
\newcommand{\rell}[3]{ s_{#1} = s_{#2} s_{#3} s_{#2}^{-1}}
\newcommand{\relli}[3]{ s_{#1} = s_{#2}^{-1} s_{#3} s_{#2}}

\begin{gather}
\rell 1{10}2, \quad
\relli 2{9}3,\quad
\relli 3{6}4,\quad
\relli 4{7}5,\quad \lb{f:relations1}\\
\rell 5{11}{6}, \quad
\relli 6{4}7,\quad
\relli 7{2}8,\quad
\rell 8{11}9, \quad \lb{f:relations2}\\
\relli 9{7}{10},\quad
\rell {10}{8}{11}, \quad
\rell {11}{5}1.\lb{f:relations3}
\end{gather}

There is a \ho~ $h:\pi_1(\gC)\to S(5)$
given by the following formulas:
\begin{gather}
h(s_1)=h(s_5)=h(s_6)=h(s_{11})=(253),\quad
h(s_2)=(234),\quad \\
h(s_3)=h(s_7)=(123),\quad
h(s_4)=(135),~
h(s_8)=(142),\\
h(s_9)=(145),~
h(s_{10})=(345).
\end{gather}
The image of $h$ is contained in the subgroup $A(5)$.
The group $S(5)$ acts by permutation of coordinates
on the free $\ZZZ$-module $\ZZZ^5$
and we obtain therefore a representation $\rho:
\pi_1\gC\to SL(n, \ZZZ)$.
The twisted Novikov homology $H_1(C_\gK; \rho^{*})$, 
with respect to the conjugate representation $\r^{*}$, 
can be computed from the free $\ZZZ((t))$-chain complex
$$C_0\lTo^{\pr_1} C_1\lTo^{\pr_2} C_2$$
where $\rk C_0=5, \rk C_1=55=\rk C_2$.
The generators of $C_1$ correspond to $s_i, 1\leq i \leq 11$,
the generators of $C_2$ correspond to
the eleven relations 
\rrf{f:relations1} - \rrf{f:relations3}, and
the matrix of $\pr_2$ is obtained by the Fox
calculus using these relations.
The homomorphism $\pr_1$ is epimorphic,
since the 0-dimensional twisted
Novikov homology of a knot 
always vanishes. Thus the rank of the $55\times 55$-matrix
of $\pr_2$ is not more than 50. The determinant of the
$50\times 50$-minor of the matrix of $\pr_2$ 
%can be computed
%with the help of a computer, it is equal to
obtained from the matrix by omitting the last 
five columns and the last five rows is equal to

\begin{gather*}-5t^{-29}+14t^{-28}-15t^{-27}+16t^{-26}
-19t^{-25}+10t^{-24}+5t^{-23}-24t^{-22}\\
+34t^{-21}-32t^{-20}+34t^{-19}\\
-24t^{-18}+5t^{-17}+10t^{-16}-19t^{-15}
+16t^{-14}-15t^{-13}+ 14t^{-12}-5t^{-11}.
\end{gather*}

This polynomial is a non-invertible element of $\ZZZ((t))$,
since the leading  coefficient is $-5\not=\pm 1$.
Therefore the torsion part of the twisted Novikov homology
in dimension $1$ is not zero, and
$q_1(\gC; \r)>0$. By Corollary
\mrf{c:n_times}
we deduce the inequality
\rrf{f:conw_estimate}.
$\qs$

By using the Kodama's KNOT program, 
we can show that the Kinoshita-Terasaka knot $\gK\gT$
has also a representation 
$\rho:
\pi_1\gK\gT\to SL(5, \ZZZ)$ 
such that $\wh q_1(\gK\gT,\r)\not= 0$.

%%%%%%%%%%%%%%%%%%%%%%%%%%%%%%%%%%%%%%%%%%%%%%%%%%%%

\section{The Kinoshita-Terasaka and Conway knots:
upper bounds for the Morse-Novikov numbers}
\mlb{s:kt_c2}

In this section, we show that both  
$\MM\NN(\gK\gT)$ and $\MM\NN(\gC)$ are 
less than or equal to 2.
Therefore 
$\MM\NN(n\gK\gT)=\MM\NN(n\gC)\leq 2n$ 
by the inequality \rrf{f:connect}.

Here we use the minimal genus Seifert surfaces 
for $\gK\gT$ and $\gC$, which were found in \cite{gabai2}.  
See Figures 3 and 4. 
Since the proofs are same, we consider only the Conway knot 
$\gC$.

\begin{figure}
\centering
\includegraphics[width=.6\textwidth]{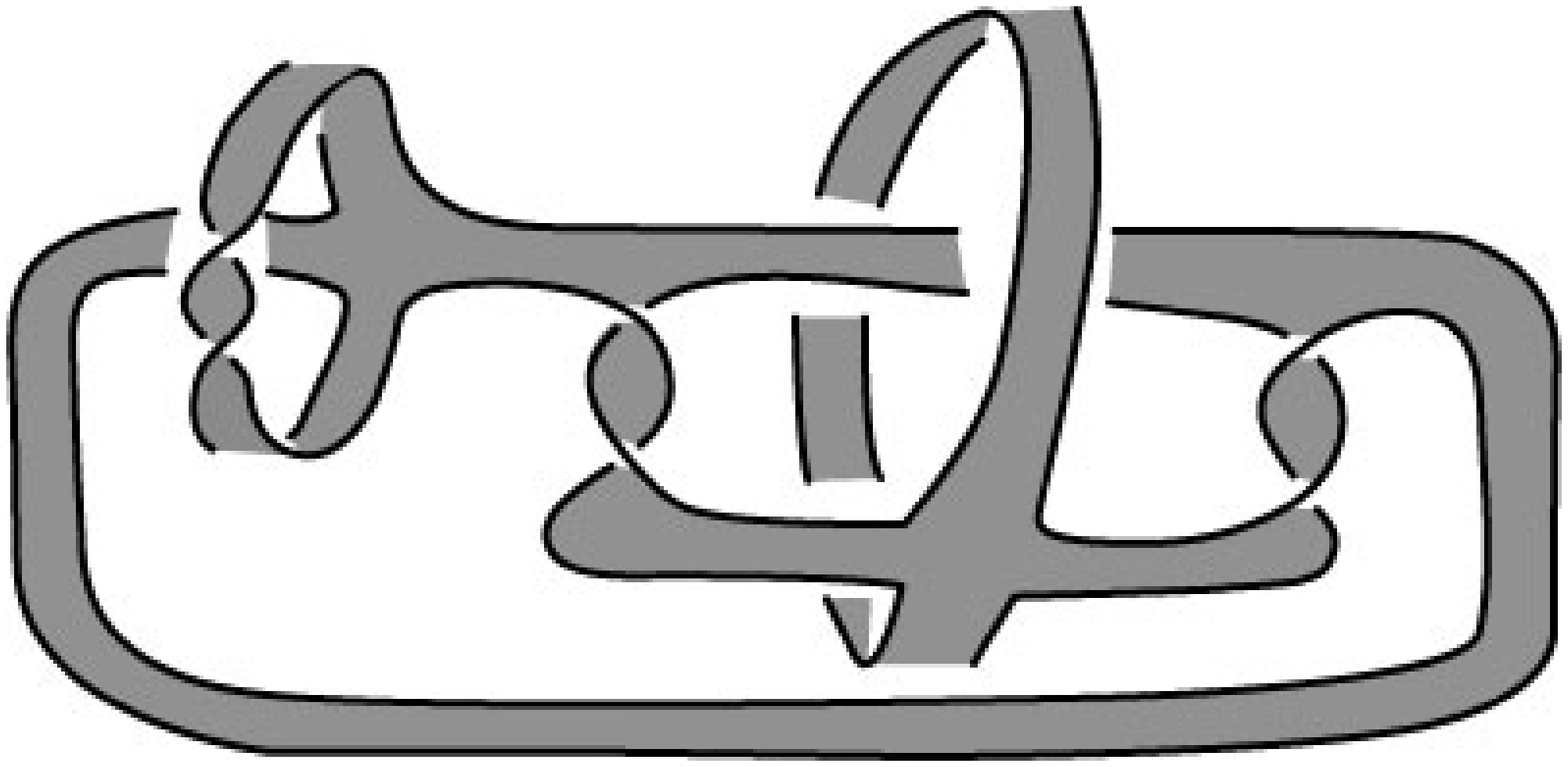}
\caption{The Kinoshita-Terasaka knot $\gK\gT$}
\end{figure}

\begin{figure}
\centering
\includegraphics[width=.6\textwidth]{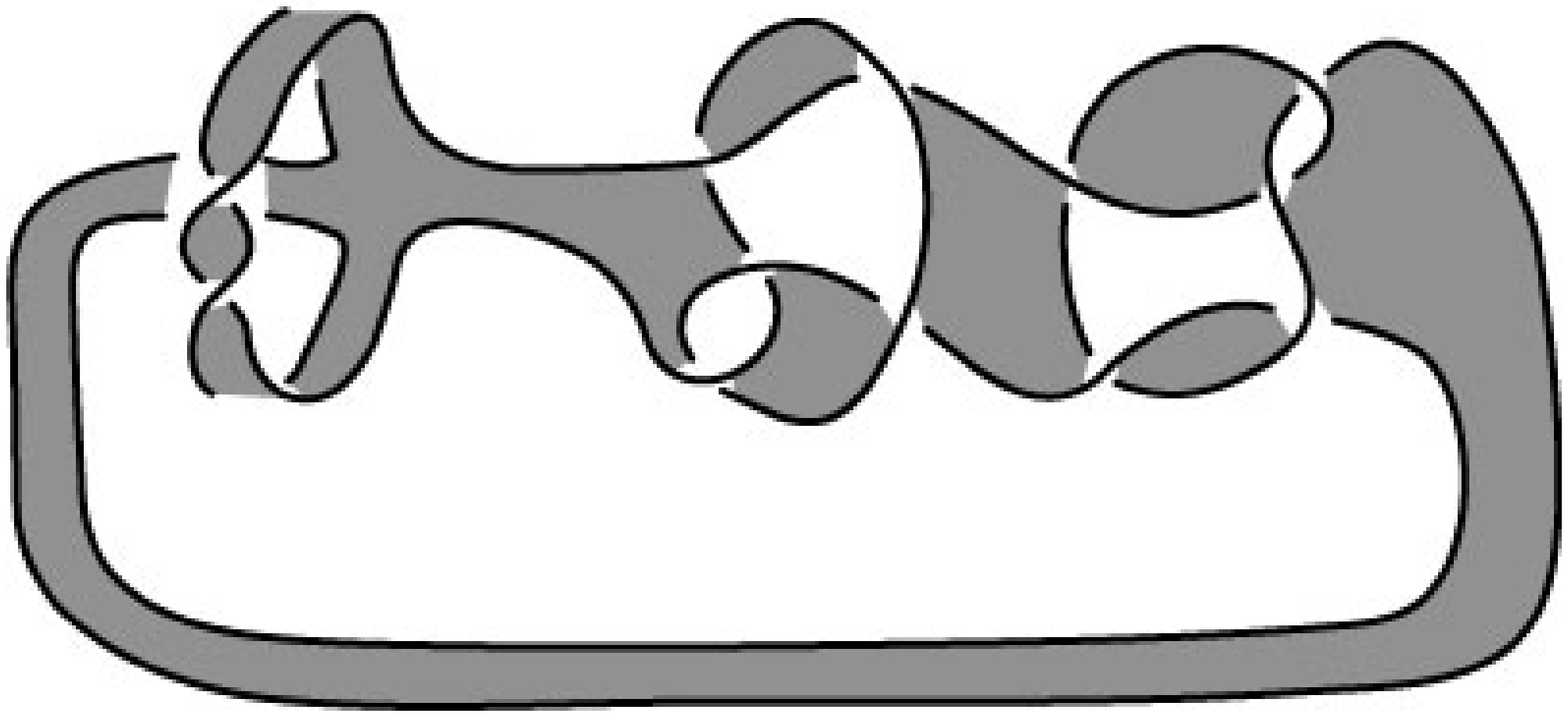}
\caption{The Conway knot $\gC$}
\end{figure}

\begin{figure}
\centering
\includegraphics[width=.7\textwidth]{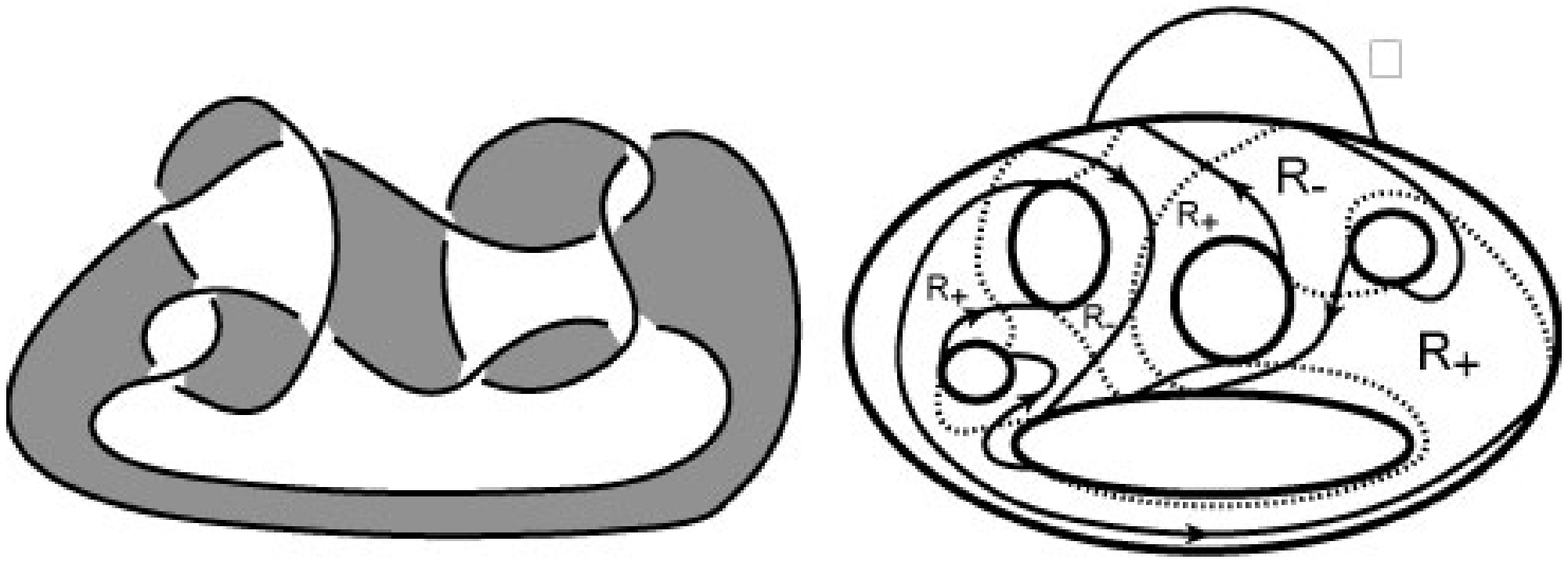}
\caption{}
\end{figure}

Since the Hopf band is a fiber surface, we may calculate 
the handle number of the Seifert surface illustrated 
in left-hand Figure 5 by Theorem B in \cite{godatwo}. 
We call $R$ this Seifert surface, and denote by $(M,\gamma)$ 
the sutured manifold for $R$. 
Further, we name $R_{+}(\g)$ and $R_{-}(\g)$, and 
let $\alpha$ be an arc properly embedded in $M$ 
as in Figure 5. (We abbreviate $R_{\pm}(\g)$ to 
$R_{\pm}$, and $\g$ to its core circles.)
Then, $N(R_{+}(\g)\cup\alpha, M)$ becomes a compression body $W$ 
with $h(W)=1$. 
By using Lemma 2.4 in \cite{godaone}, we can observe 
that cl$(M-W)$ is a compression body $W'$
such that $\partial_{-}W'=R_{-}(\g)$ and $h(W')=1$. 
Hence we have $h(R)\le 1$, namely, 
$\MM\NN(\gC)\le 2$. 
This completes the proof.

%%%%%%%%%%%%%%%%%%%%%%%%%%%%%%%%%%%%%%%%%%%%%%%%%%%%%%%%
%\input tw_alex.doc

\section{Relation with the twisted Alexander invariants}
\mlb{s:tw_alex}

In this section we investigate the relations between
the invariants arising from the Morse-Novikov theory and
the twisted Alexander invariant.
We shall use the version of 
twisted Alexander polynomials 
via the Reidemeister torsion 
as suggested by Kitano in \cite{kitano}. 
It is proved in \cite{kitano} 
that this version is equivalent to the 
initial Wada's definition. 
Let $K$ be an oriented knot in $S^3$.
Let $\xi:G=\pi_1K\to \ZZZ$
the homomorphism introduced in the previous
section; let $\r:G\to GL(n,F)$
be a right representation of the knot group, where 
$F$ is a field. Let $S=F(t)$;
the representation $\r_\xi:G\to SL(n,Q)$
determines also the representation
$\wi\r_\xi:G\to GL(n,S)$.
Form the chain complex
\bq\lb{f:chain_rat}
C_*(K;\r )=
S^n\tens{\wi\r_\xi} ~C_*(\wi{C}_K).
\end{equation}
This is a chain complex of vector spaces over $S$
and the same is true for its homology.
We shall limit ourselves to the case when
the chain complex \rrf{f:chain_rat}
is acyclic. 
In this case the {\it twisted Alexander invariant} 
is defined as the Reidemeister torsion of the chain
complex $C_*(K;\r)$.
It is an element of $S$, well-defined modulo multiplication by $\pm t^n, n\in\ZZZ$.
The field $S$ can be considered as 
a subfield of the  Novikov ring $F((t))$, and therefore 
the twisted Alexander invariant
can be also considered as an element of 
$F((t))$, that is, as  a  Laurent series,
well-defined up to multiplication by $\pm t^n$.
This Laurent series will be denoted by 
$\tau(K; \rho)$. 

\bede\mlb{d:monic}
An element $x\in F((t))$
will be called {\it monic}, if the coefficient 
corresponding to the term of the lowest degree of $t$ equals $1$ or $-1$.
\end{defi}

\bepr\mlb{p:nov_alex}
Assume that 
the representation $\rho:\pi_1K\to SL(n,F)$
is obtained from a representation
$\rho_0:\pi_1K\to SL(n,\ZZZ)$
via the tautological homomorphism
$\ZZZ\to F, ~1\to 1$.
Assume that
the Novikov homology  $\wh H_*(K, \rho_0)$
vanishes.
Then the twisted Alexander polynomial $\tau(K,\rho)$
is $\pm$-monic.
\enpr
\Prf
The chain complex $\wh C_*(C_K, \rho_0)$
is acyclic. The Whitehead   torsion of this complex
is an element $$
\tau_0\in \ove K_1\big(~\wh \ZZZ((t))~\big)\approx \ZZZ((t))^\bu,
$$
defined up to multiplication by $\pm t^n$.
The multiplicative group  $\ZZZ((t))^\bu$
of all invertible  elements of $\ZZZ((t))$
consists of the Laurent series $\l$ such that 
the lowest degree term equals $\pm t^n$.
The image of  $\tau_0$ in
$F((t))$ coincides with the image
of the twisted Alexander polynomial
$\tau(K, \rho)$, therefore the latter is monic.
$\qs$

%\bigskip
%Example.  

%%%%%%%%%%%%%%%%%%%%%%%%%%%%%%%%%%%%%%%%%%%%%%%%

%\input conj.doc

\section{A conjecture}
\mlb{s:conj_nov}

The examples of the previous sections
show that the lower bounds provided by the
twisted Novikov
homology are often optimal.
Every such estimate arises from  a
linear representation
of a fundamental group.
In the present section we  explore
an approach which starts from the
most general form of the Novikov homology
and which should give the best possible
lower bounds.

Let us recall first the construction of
the Novikov completion of a given group $G$ \wrt~
a \ho~ $\xi:G\to \ZZZ$.
Set $\L=\ZZZ G$
and denote by
$\wh{\wh\L}$
the abelian group of all functions $G\to\ZZZ$.
Equivalently,
$\wh{\wh\L}$
is the set of all formal linear combinations
$\l=\sum_{g\in G} n_g g$
 (not necessarily finite)
of the elements of $G$ with integral coefficients.
For $\l\in \wh{\wh\L}$ set
$\supp\l=
\{g\in G\mid n_g\not= 0\}$.
Set
\bq
\wh\L_{\xi}=
\{\l\in
\wh{\wh \L}
\mid \forall C\in\RRR \quad
\supp\l\cap\xi^{-1}([C,\infty[)
\mbox{ is finite }\}
\end{equation}
Then $\Lxi$ has a natural structure of  a ring.

The following theorem was
first proved in \cite{patou}.

\beth\mlb{t:nov_fibr}$($\cite{patou}$)$
Let $f:M\to S^1$ be a Morse map
such that $f_*=\xi:\pi_1M\to\ZZZ$.
Then there is a chain complex $\NN_*$
of  free finitely generated
$\Lxi$-modules (the Novikov complex), such that
\been\item
 the number of
free generators in each degree $i$
equals $m_i(f)$,
\item
there is a chain homotopy equivalence
$$
\phi: \NN_*
\rTo
\Lxi\tens{\L}C_*^\D(\wi M).
$$
\enen
\enth
This theorem was proved 
in \cite{patou}
for the case of closed manifolds.
A similar result holds for the regular Morse maps
$f:C_L\to S^1$.

\begin{rema}
In the paper \cite{patou},
we worked with the convention 
that the fundamental group 
acts on the universal covering 
on the right. 
Theorem \mrf{t:nov_fibr} above 
is the translation of the results 
of \cite{patou} to the 
language of the left modules. 
\end{rema}

\beco\mlb{c:nov_fibr}
If the link $L$ is fibred, then
$$
H_*\big(\Lxi\tens{\L}C_*^\D(\wi C_L)\big)=0.
$$
\enco

In general the vanishing of the Novikov homology
is not sufficient for existence of a fibration.
Even in the case when the dimension
of the manifold is $\geq 6$
there is a secondary obstruction
in the Whitehead group of the
Novikov ring of the fundamental group.
Let us give the definition.
An element $\l\in\Lxi$
is called {\it trivial unit } if
$$\l= \ve\cdot g\cdot (1 +\mu)
\mxx{ with }
\ve=\pm 1, ~ g\in G,~ \supp \mu\in
\xi^{-1}(]-\infty,0[). $$
The set of all trivial units
is a multiplicative subgroup $U$
of the group of all units of $\Lxi$.
The group $\ove K_1(\Lxi)/U'$
where $U'$ if the image of $U$, is
denoted $\Wh(G,\xi)$ and called {\it the Whitehead
group of the  Novikov ring}.
The next theorem follows immediately
by combining the classical theorem of Waldhausen
\cite{waldhausen}
(the Whitehead group of the
link group vanishes) and the main theorem
of \cite{pajandran}.

\beth\mlb{t:wh_knots_vanish}
We have $\Wh(\pi_1L, \xi)=0$
where $\xi:\pi_1L\to\ZZZ$
is the homomorphism which sends
each meridian generator to $1\in \ZZZ$. $\qs$
\enth

Thus in the case of links
the total obstruction to fibering provided
by the Novikov complex is reduced
to the Novikov homology.

\begin{conj}\mlb{conj:fibr}
A link $L$ is fibred if and only if
\bq\lb{f:cond_fib_nov}
H_*\big(\Lxi\tens{\L}C_*^\D(\wi C_L)\big)
=0.
\end{equation}
\end{conj}

The Novikov ring $\Lxi$
is a complicated algebraic object, and
the verification of the condition
\rrf{f:cond_fib_nov}
is certainly a difficult algebraic task.
The twisted Novikov homology as introduced and studied
in the previous sections
provides an effectively computable tools
for evaluating the Novikov homology, and
as we have seen in many examples,
the twisted Novikov homology is often sufficient
to compute the Morse-Novikov number.
%Thus our next conjecture is the following.
Thus we are led to the following problem.

\begin{prob}\mlb{prob:second}
Is it true that 
vanishing of the $\r$-twisted Novikov homology 
for every representation implies 
the condition \rrf{f:cond_fib_nov} ?
%If the $\rho$-twisted Novikov homology vanishes
%for every representation $\rho:\pi_1L\to GL(n,R)$
%with $R=\ZZZ$ or $R$ is a field,
%then we have
%\rrf{f:cond_fib_nov}.
\end{prob}

A natural and a very interesting question
would be to investigate the relations
between the Problem \mrf{prob:second}
and the Problem 1.1 of \cite{gomo},
which asks whether the information
contained in the twisted Alexander
polynomials for all $SL(2k,\FFF)$-representations
is sufficient to decide whether a link is fibred.
%Observe that the comparison of these two
%problems may boil down to  purely algebraic
%questions.

\bibliographystyle{amsplain}

\begin{thebibliography}{30}


\bibitem{bz}
G. Burde, H. Zieschang, 
Knots,  
de Gruyter Studies in Mathematics, 5. 
Walter de Gruyter $\&$ Co., Berlin, 1985. 


\bibitem{CG}
A. Casson and  C. McA. Gordon,
\textit{Reducing Heegaard splittings},  
Topology Appl. 27 (1987), no. 3, 275--283. 


\bibitem{conway}
J. H. Conway,
\textit{
An enumeration of 
knots and links, and some 
of their algebraic properties},
Comutational 
Problems in Abstract Algebra,
Pergamon Press, 
New York, 1970,  329--358.

\bibitem{gabai} 
D. Gabai, 
\textit{Foliations and the topology of $3$-manifolds}, 
J. Differential Geom. 18 (1983), 445-503. 

\bibitem{gabai2}
D. Gabai,
\textit{Foliations and genera of links}, 
Topology 23 (1984), 381-394.

\bibitem{godaone}
H. Goda,
\textit{Heegaard splitting for sutured manifolds and Murasugi sum},
Osaka J. Math. 29 (1992), 21-40.

\bibitem{godatwo}
H. Goda,
\textit{On handle number of Seifert surfaces in $S^{3}$}, 
Osaka J. Math. 30 (1993), 63-80.


\bibitem{gomo}
H. Goda,  T. Morifuji,
\textit{Twisted Alexander polynomial
for $SL(2,\CCC)$-representations
and fibered knots},
C. R. Math. Acad. Sci. Soc. R. Canada 25, (2003), 97-101. 


\bibitem{gokimo}
H. Goda, T. Kitano, T. Morifuji,
\textit{Reidemeister torsion,
twisted Alexander polynomial
and fibered knots},
preprint, 2002.

\bibitem{hirarud}
M. Hirasawa, L. Rudolph, 
\textit{Constructions of Morse maps for knots ans 
links, and upper bounds on the Morse-Novikov number}, 
preprint, 2003. 

\bibitem{kinoshitaterasaka} 
S. Kinoshita, H. Terasaka, 
\textit{On unions of knots}, 
Osaka Math. J. 9 (1957), 131-153. 

\bibitem{kirklivingston}
P. Kirk, C. Livingston, 
\textit{Twisted Alexander invariants, Reidemeister torsion, 
and Casson-Gordon invariants}, 
Topology 38 (1999), 635-661.

\bibitem{kitano}
T. Kitano,
\textit{
Twisted Alexander polynomial
and  Reidemeister torsion}, 
Pacific J.Math,
174 (1996), 431-442.


\bibitem{lin}
X.S. Lin, 
\textit{Representations of knot groups and twisted Alexander polynomials}, 
Acta Math. Sin. (Engl. Ser.) 17 (2001), 361--380.

\bibitem{novidok}
S. P. Novikov, 
\textit{Multivalued functions and functionals, 
An analogue of the Morse theory}, 
Soviet Math.Dokl.24 (1981), 222-226. 

%\bibitem{noviuspe}
%S. P. Novikov,
%\textit{The hamiltonian formalism and a
%multivalued analogue of
%Morse theory},  
%Russ. Math. Surveys 37 (1982),
%1-56.

\bibitem{patou}
A.V. Pajitnov, 
\textit{On the Novikov
complex for rational Morse forms},
preprint: Institut for Mathematik og datalogi, 
Odense Universitet Preprints 1991, 
No 12, Oct. 1991; 
http://193.52.98.6/~Pajitnov/od.pdf
journal article: 
Annales de la Facult\'e de Sciences de
Toulouse 4 (1995), 297-338.

\bibitem{pajandran}
A. V. Pajitnov, A. Ranicki,
\textit{The Whitehead group of the Novikov ring}, 
K-theory 21 (2000), 325--365. 

\bibitem{prw}
A. V. Pajitnov, L. Rudolph, C. Weber, 
\textit{Morse-Novikov number for knots and links}, 
St. Petersburg Math. J. 13 (2002), 417-426.

\bibitem{rolfsen}
D. Rolfsen, 
Knots and links,
Mathematics Lecture Series, No. 7. 
Publish or Perish, Inc., Berkeley, Calif., 1976.


\bibitem{wada} M. Wada,
\textit{Twisted Alexander polynomial for finitely 
presentable groups},
Topology 33 (1994), 241--256.


\bibitem{waldhausen} F. Waldhausen,
\textit{Algebraic $K$-theory of generalized free products. I, II},
Ann. of Math. (2) 108 (1978),  135--204.

\end{thebibliography}

\end{document}